\documentclass[12pt, reqno]{amsart}
\usepackage{amssymb, amsthm, mathrsfs}

\newtheorem{thm}{Theorem}

\newtheorem{lemma}[thm]{Lemma}
\newtheorem{theorem}[thm]{Theorem}
\newtheorem{prop}[thm]{Proposition}

\numberwithin{equation}{section}

\theoremstyle{definition}

\newcommand{\al}{\alpha}
\renewcommand{\b}{\beta}

\newcommand{\e}{\varepsilon}
\newcommand{\la}{\lambda}
\renewcommand{\phi}{\varphi}

\renewcommand{\d}{\partial}
\newcommand{\R}{{\mathbb R}}

\newcommand{\Case}[1]{\noindent \textbf{Case #1.}}

\newcommand{\nequiv}{\mathrel{\setbox0\hbox{$\equiv$}%
                     \rlap{\hbox{$\equiv$}}\hbox to \wd0{\hfil $/$\hfil}}}

\renewcommand{\qed}{\rule{3mm}{3mm}}
\renewenvironment{proof}
    {\vspace{1mm}\noindent\textbf{Proof.}}
    {\hspace*{\fill} $\qed$\vspace{1mm}}

\begin{document}
\title[Properties of solutions to the nonlinear biharmonic equation]
    {Stability and intersection properties of solutions to the nonlinear biharmonic equation}
\author{Paschalis Karageorgis}
\address{School of Mathematics, Trinity College, Dublin 2, Ireland.}
\email{pete@maths.tcd.ie}

\begin{abstract}
We study the positive, regular, radially symmetric solutions to the nonlinear biharmonic equation $\Delta^2 \phi =
\phi^p$. First, we show that there exists a critical value $p_c$, depending on the space dimension, such that the
solutions are linearly unstable if $p<p_c$ and linearly stable if $p\geq p_c$.  Then, we focus on the supercritical case
$p\geq p_c$ and we show that the graphs of no two solutions intersect one another.
\end{abstract}
\maketitle

\section{Introduction}
Consider the positive, regular, radially symmetric solutions of the equation
\begin{equation}\label{ee}
\Delta^2 \phi (x) = \phi(x)^p, \quad\quad x\in \R^n.
\end{equation}
Such solutions are known to exist when $n>4$ and $p\geq \frac{n+4}{n-4}$, but they fail to exist, otherwise. Our main
goal in this paper is to study their qualitative properties, and to also relate those to the well-understood properties
of solutions to the second-order analogue
\begin{equation}\label{2ee}
-\Delta \phi (x) = \phi(x)^p, \quad\quad x\in \R^n.
\end{equation}

Linear stability for the second-order equation \eqref{2ee} was addressed in \cite{KaSt}, where the positive, regular,
radially symmetric solutions were found to be linearly stable if and only if
\begin{equation}\label{sc2}
p \cdot Q_2\left( \frac{2}{p-1} \right) \leq Q_2\left( \frac{n-2}{2} \right),
\quad\quad Q_2(\al) \equiv |x|^{\al+2} (-\Delta) \,|x|^{-\al}.
\end{equation}
In this paper, we establish a similar result for the fourth-order equation \eqref{ee}, namely that the positive, regular,
radially symmetric solutions are linearly stable if and only if
\begin{equation}\label{sc}
p \cdot Q_4\left( \frac{4}{p-1} \right) \leq Q_4\left( \frac{n-4}{2} \right),
\quad\quad Q_4(\al) \equiv |x|^{\al+4} \,\Delta^2 |x|^{-\al}.
\end{equation}
Although originally stated in a different way, the last two conditions appeared in the work of Wang~\cite{Wang} for the
second-order equation and Gazzola and Grunau \cite{GG} for the fourth-order one.  Among other things, these authors
studied the intersection properties\footnote{Due to scaling, the existence of one solution implies the existence of
infinitely many solutions.} of radially symmetric solutions, and they found that the above conditions play a crucial role
in that context.

According to a result of Wang \cite{Wang} for the second-order equation, the graphs of no two radially symmetric
solutions intersect one another, if \eqref{sc2} holds, while the graphs of any two radially symmetric solutions intersect
one another, otherwise.  Although a similar dichotomy is expected to hold for the fourth-order equation, we are only able
to prove the first part of such a result, namely that the graphs of no two radially symmetric solutions of \eqref{ee}
intersect one another, if \eqref{sc} holds.  This already improves a result of \cite{GG}, which shows that the number of
intersections is at most finite, if \eqref{sc} holds with strict inequality.  As for the remaining case in which
\eqref{sc} is violated, neither our approach nor the one in \cite{GG} provides any conclusions.

In section~\ref{li}, we show that \eqref{sc} is a necessary condition for linear stability.  In section~\ref{ls}, we show
that it is also a sufficient condition.  Our main results appear in section \ref{omr}, where we also simplify the
stability condition \eqref{sc}.  Our stability result is given in Theorem \ref{main}, and our result on the intersection
properties of solutions is given in Theorem \ref{sep}.

\section{Linear instability}\label{li}
The main result in this section is Proposition \ref{unst}, which gives a sufficient condition for linear instability.
Although this condition will be simplified in section \ref{omr}, it is much more convenient to initially state it in
terms of the quartic polynomial
\begin{equation}\label{Q}
Q_4(\al) \equiv |x|^{\al+4} \,\Delta^2 |x|^{-\al} = \al (\al+2) (\al+2-n) (\al+4-n).
\end{equation}
This polynomial is closely related to Rellich's inequality
\begin{equation}\label{Rel}
\int_{\R^n} (\Delta u)^2 \:dx \geq \frac{n^2(n-4)^2}{16} \:\int_{\R^n} |x|^{-4} u^2\:dx,
\end{equation}
which is valid for each $u\in H^2(\R^n)$ and each $n>4$.  Namely, the constant that appears on the right hand side is
merely the unique local maximum value of $Q_4$, and it is known to be sharp in the following sense.

\begin{lemma}\label{Rell}
Let $n>4$ and let $V$ be a bounded function on $\R^n$ that vanishes at infinity.  If there exists some $\e>0$ such that
\begin{equation*}
V(x) \leq -(1+\e)\cdot \frac{n^2(n-4)^2}{16} \cdot |x|^{-4}
\end{equation*}
for all large enough $|x|$, then the operator $\Delta^2+V$ has a negative eigenvalue.
\end{lemma}

For a proof of Rellich's inequality \eqref{Rel} and Lemma \ref{Rell}, we refer the reader to section~II.7 in Rellich's
book \cite{Rell}.  We now use the previous lemma to address the linear instability of positive, regular solutions to
\eqref{ee}. The known results on the existence of such solutions are summarized in our next lemma; see \cite{WX,CLO,GG}
for parts (a), (b) and (c), respectively.

\begin{lemma}\label{ex}
Let $n\geq 1$ and $p>1$.  Denote by $Q_4$ the quartic in \eqref{Q}.

\begin{itemize}
\item[(a)]
If either $n\leq 4$ or $p<\frac{n+4}{n-4}$, then equation \eqref{ee} has no positive $\mathcal{C}^4$ solutions.

\item[(b)]
If $n>4$ and $p= \frac{n+4}{n-4}$, then all positive $\mathcal{C}^4$ solutions of \eqref{ee} are of the form
\begin{equation}\label{phila}
\phi_\la(x) = \Bigl[ (n-4)(n-2)n(n+2) \Bigr]^{\frac{1}{p-1}} \cdot
\left( \frac{\la}{\la^2 + |x-y|^2} \right)^{\frac{4}{p-1}}
\end{equation}
for some $\la>0$ and some $y\in\R^n$.

\item[(c)]
If $n>4$ and $p>\frac{n+4}{n-4}$, then the positive $\mathcal{C}^4$, radially symmetric solutions of \eqref{ee} form an
one-parameter family $\{ \phi_\al \}_{\al>0}$, where each $\phi_\al$ is such that
\begin{equation}\label{phia}
\phi_\al(0)= \al, \quad\quad \lim_{|x|\to \infty} |x|^4 \,\phi_\al(x)^{p-1} = Q_4\left( \frac{4}{p-1} \right) > 0.
\end{equation}
\end{itemize}
\end{lemma}

\begin{prop}\label{unst}
Let $n>4$ and $p\geq \frac{n+4}{n-4}$.  Let $Q_4$ be the quartic in \eqref{Q} and let $\phi$ denote any one of the
solutions provided by Lemma \ref{ex}.  Then $\Delta^2 - p\phi^{p-1}$ has a negative eigenvalue, if
\begin{equation}\label{ic}
p \cdot Q_4\left( \frac{4}{p-1} \right) > Q_4\left( \frac{n-4}{2} \right) = \frac{n^2(n-4)^2}{16} \,.
\end{equation}
In particular, it has a negative eigenvalue, if $p=\frac{n+4}{n-4}$.
\end{prop}

\begin{proof}
Suppose first that $p>\frac{n+4}{n-4}$.  Using part (c) of Lemma \ref{ex} and our assumption \eqref{ic}, we can then find
some small enough $\e>0$ such that
\begin{equation*}
\lim_{|x|\to \infty} |x|^4 \,\phi(x)^{p-1} = Q_4\left( \frac{4}{p-1} \right)> p^{-1}(1+2\e)\cdot\frac{n^2(n-4)^2}{16} \,.
\end{equation*}
Since this implies that
\begin{equation*}
V(x) \equiv -p\phi(x)^{p-1} < -(1+\e)\cdot \frac{n^2(n-4)^2}{16} \cdot |x|^{-4}
\end{equation*}
for all large enough $|x|$, the existence of a negative eigenvalue follows by Lemma \ref{Rell}.

Suppose now that $p=\frac{n+4}{n-4}$.  Then our assumption \eqref{ic} automatically holds because
\begin{equation*}
\frac{n^2(n-4)^2}{16} = Q_4\left( \frac{n-4}{2} \right) = Q_4\left( \frac{4}{p-1} \right) < p\cdot Q_4\left(
\frac{4}{p-1} \right)
\end{equation*}
for this particular case.  According to part (b) of Lemma \ref{ex}, we also have
\begin{equation*}
\Delta^2 - p\phi(x)^{p-1}= \Delta^2 -\la^4 (n-2)n(n+2)(n+4) \cdot \left( \la^2 + |x-y|^2 \right)^{-4}
\end{equation*}
for some $\la>0$ and some $y\in \R^n$.  Thus, it suffices to check that the associated energy
\begin{equation*}
E(\zeta) = \int_{\R^n} (\Delta \zeta)^2 \:dx - \int_{\R^n} p\phi^{p-1} \zeta^2 \:dx
\end{equation*}
is negative for some test function $\zeta\in H^2(\R^n)$.  Let us then consider the test function
\begin{equation*}
\zeta(x) = \left( \la^2 + |x-y|^2 \right)^{-\frac{n-2}{2}}.
\end{equation*}
Since $n>4$, we have $\zeta\in H^2(\R^n)$, while a straightforward computation gives
\begin{align*}
E(\zeta) &= -8\la^4 n(n-2)(n+1) \int_{\R^n} \frac{dx}{(\la^2 +|x-y|^2)^{n+2}} <0.
\end{align*}
This implies the presence of a negative eigenvalue and it also completes the proof.
\end{proof}

\section{Linear stability}\label{ls}
In this section, we address the linear stability of the solutions provided by Lemma \ref{ex}.  First, we use an
Emden-Fowler transformation to transform \eqref{ee} into an ODE whose linear part has constant coefficients. Although
this transformation is quite standard, the subsequent part of our analysis is not. The main result of this section is
given in Proposition \ref{st}.

\begin{lemma}\label{cco}
Let $n\geq 1$ and $p>1$.  Let $Q_4$ be the quartic in \eqref{Q} and suppose $\phi$ is a positive solution of the
biharmonic equation \eqref{ee}. Setting $m=\frac{4}{p-1}$ for convenience, the function
\begin{equation}\label{W}
W(s) = e^{m s} \phi(e^s) = r^m \phi(r), \quad\quad s= \log r= \log |x|
\end{equation}
must then be a solution to the ordinary differential equation
\begin{equation}\label{ne1}
Q_4(m-\d_s) \,W(s) = W(s)^p.
\end{equation}
\end{lemma}

\begin{proof}
Since $\d_r = e^{-s}\d_s$, a short computation allows us to write the radial Laplacian as
\begin{equation*}
\Delta = \d_r^2 + (n-1) r^{-1}\d_r = e^{-2s} (n-2+\d_s) \d_s.
\end{equation*}
Using the operator identity $\d_s e^{-ks} = e^{-ks}(\d_s-k)$, one can then easily check that
\begin{align*}
\Delta^2 e^{-m s} &= e^{-4s -m s} \,Q_4(m-\d_s) = e^{-m p s} \,Q_4(m-\d_s).
\end{align*}
This also implies that $Q_4(m-\d_s) \,W(s) = e^{m ps} \Delta^2 \phi (e^s) = W(s)^p$, as needed.
\end{proof}

\begin{lemma}\label{poly}
Let $n>4$ and $p>\frac{n+4}{n-4}$.  Set $m=\frac{4}{p-1}$ and let $Q_4$ be the quartic in \eqref{Q}.  Assuming that the
stability condition \eqref{sc} holds, the polynomial
\begin{equation}\label{ch1}
\mathscr{P}(\la) = Q_4(m-\la) -pQ_4(m)
\end{equation}
must then have four real roots $\la_1,\la_2,\la_3<0<\la_4$.
\end{lemma}

\begin{proof}
Noting that $Q_4$ is symmetric about $\frac{n-4}{2}$, we see that $\mathscr{P}$ is symmetric about
\begin{equation*}
\la_* \equiv m- \frac{n-4}{2} = \frac{4}{p-1} - \frac{n-4}{2} \,,
\end{equation*}
where $\la_* < 0$ because $p>\frac{n+4}{n-4}$ by assumption.  Moreover, we have
\begin{equation*}
\lim_{\la\to \pm\infty} \mathscr{P}(\la) = + \infty, \quad\quad
\mathscr{P}(2\la_*) = \mathscr{P}(0) = (1-p) \cdot Q_4(m) < 0
\end{equation*}
because of \eqref{phia}, and we also have
\begin{equation*}
\mathscr{P}(\la_*) = Q_4\left( \frac{n-4}{2} \right) - p\cdot Q_4\left( \frac{4}{p-1} \right) \geq 0
\end{equation*}
because of \eqref{sc}.  This forces $\mathscr{P}(\la)$ to have at least one root in each of the intervals
\begin{equation*}
(-\infty,2\la_*), \quad\quad (2\la_*,\la_*], \quad\quad [\la_*,0), \quad\quad (0,\infty).
\end{equation*}
In the case that $\la_*$ itself happens to be a root, then it must be a double root by symmetry.  In any case then,
$\mathscr{P}(\la)$ has three negative roots and one positive root, as needed.
\end{proof}

\begin{prop}\label{st}
Let $n>4$ and $p>\frac{n+4}{n-4}$.  Let $Q_4$ be the quartic in \eqref{Q} and let $\phi$ denote any one of the solutions
provided by Lemma \ref{ex}.  Assuming that \eqref{sc} holds, one has
\begin{equation}\label{bel1}
|x|^4 \,\phi(x)^{p-1} \leq Q_4\left( \frac{4}{p-1} \right)
\end{equation}
for each $x\in\R^n$, and the operator $\Delta^2-p\phi^{p-1}$ has no negative spectrum.
\end{prop}

\begin{proof}
First, suppose that \eqref{bel1} does hold.  Using our assumption \eqref{sc}, we then get
\begin{equation*}
-p\phi(x)^{p-1} \geq -p \cdot Q_4\left( \frac{4}{p-1} \right) \cdot |x|^{-4} \geq -\frac{n^2(n-4)^2}{16} \cdot |x|^{-4}
\end{equation*}
for each $x\in\R^n$, so $\Delta^2-p\phi^{p-1}$ has no negative spectrum by Rellich's inequality \eqref{Rel}.

Let us now focus on the derivation of \eqref{bel1}.  Set $m=\frac{4}{p-1}$ and consider the function
\begin{equation*}
W(s) = e^{ms} \phi(e^s)= r^m \phi(r), \quad\quad s=\log r= \log |x|.
\end{equation*}
Then $W(s)$ is positive and it satisfies the equation
\begin{equation}\label{11}
Q_4(m-\d_s) \,W(s) = W(s)^p
\end{equation}
by Lemma \ref{cco}. We note that $s$ ranges over $(-\infty,\infty)$ as $r$ ranges from $0$ to $\infty$, while
\begin{equation*}
\lim_{s\to -\infty} W(s) = \lim_{r\to 0^+} r^m\phi(r) = 0.
\end{equation*}
The derivatives of $W(s)$ must also vanish at $s=-\infty$ because
\begin{equation*}
\lim_{s\to -\infty} W'(s) = \lim_{r\to 0^+} r\cdot \d_r [r^m\phi(r)] = 0,
\end{equation*}
and so on. Using the fact that $x\mapsto x^p$ is convex on $(0,\infty)$, we now find
\begin{align}\label{12}
W(s)^p - Q_4(m)^{\frac{p}{p-1}} &\geq pQ_4(m)\cdot \left( W(s)-Q_4(m)^{\frac{1}{p-1}} \right).
\end{align}
Inserting this inequality in \eqref{11}, we thus find
\begin{equation}\label{13}
Q_4(m-\d_s) \,W(s) - Q_4(m)^{\frac{p}{p-1}} \geq pQ_4(m)\cdot \left( W(s)-Q_4(m)^{\frac{1}{p-1}} \right).
\end{equation}
To eliminate the constant term on the left hand side, we change variables by
\begin{equation}\label{Y}
Y(s) = W(s) - Q_4(m)^{\frac{1}{p-1}}.
\end{equation}
Then we can write equation \eqref{13} in the equivalent form
\begin{equation*}
\Bigl[ Q_4(m-\d_s) -pQ_4(m) \Bigr] \,Y(s) \geq 0.
\end{equation*}
Invoking Lemma \ref{poly}, we now factor the last ODE to obtain
\begin{equation}\label{14}
(\d_s - \la_1)(\d_s - \la_2)(\d_s - \la_3)(\d_s - \la_4) \,Y(s) \geq 0
\end{equation}
for some $\la_1,\la_2,\la_3<0<\la_4$.  Multiplying by $e^{-\la_1s}$ and integrating over $(-\infty,s)$, we get
\begin{equation*}
e^{-\la_1s} (\d_s - \la_2)(\d_s - \la_3)(\d_s - \la_4) \,Y(s) \geq 0
\end{equation*}
because $\la_1<0$.  We ignore the exponential factor and use the same argument twice to get
\begin{equation*}
(\d_s - \la_4) \,Y(s) \geq 0
\end{equation*}
since $\la_2,\la_3<0$ as well.  Multiplying by $e^{-\la_4s}$ and integrating over $(s,+\infty)$, we then find
\begin{equation*}
e^{-\la_4s} \Bigl[ W(s) - Q_4(m)^{\frac{1}{p-1}} \Bigr]
\leq \lim_{s\to\infty} e^{-\la_4s} \Bigl[ W(s) - Q_4(m)^{\frac{1}{p-1}} \Bigr].
\end{equation*}
The limit on the right hand side is zero because $\la_4>0$ by above and since
\begin{equation*}
\lim_{s\to \infty} W(s) = \lim_{|x|\to\infty} |x|^{\frac{4}{p-1}} \,\phi(x) = Q_4(m)^{\frac{1}{p-1}}
\end{equation*}
by \eqref{phia}. In particular, we may finally deduce the estimate
\begin{equation*}
W(s) \leq Q_4(m)^{\frac{1}{p-1}},
\end{equation*}
which is precisely the desired estimate \eqref{bel1} because $W(s) = |x|^{\frac{4}{p-1}} \phi(x)$ by above.
\end{proof}

\section{Our main results}\label{omr}
In this section, we give our main results regarding the stability and intersection properties of the positive, regular
solutions to \eqref{ee}.  Our first theorem is an easy consequence of the results obtained in the previous two sections.

\begin{theorem}\label{main}
Let $n>4$ and $p\geq p_n\equiv \frac{n+4}{n-4}$.  Let $Q_4$ be the quartic in \eqref{Q} and let $\phi$ denote any one of
the solutions provided by Lemma \ref{ex}.  Then the following dichotomy holds.

If $n\leq 12$, then $\phi$ is linearly unstable for any $p\geq p_n$ whatsoever.

If $n\geq 13$, on the other hand, then the equation
\begin{equation*}
p \cdot Q_4\left( \frac{4}{p-1} \right) = Q_4\left( \frac{n-4}{2} \right)
\end{equation*}
has a unique solution $p_c>p_n$, and $\phi$ is linearly unstable if and only if $p_c> p\geq p_n$.
\end{theorem}

\begin{proof}
Consider the expression
\begin{equation}\label{Qnew}
\mathcal{Q}(p) \equiv 16(p-1)^4\cdot \left[ Q_4\left( \frac{n-4}{2} \right) -p\cdot Q_4 \left( \frac{4}{p-1}
\right)\right].
\end{equation}
By Propositions \ref{unst} and \ref{st}, to say that $\phi$ is linearly unstable is to say that $\mathcal{Q}(p)<0$.

Let us now combine our definitions \eqref{Q} and \eqref{Qnew} to write
\begin{equation*}
\mathcal{Q}(p) = n^2(n-4)^2(p-1)^4 - 2^7p(p+1)\Bigl( (n-4)p-n \Bigr) \Bigl( (n-2)p- (n+2) \Bigr).
\end{equation*}
Using this explicit equation, it is easy to see that
\begin{equation*}
\mathcal{Q}(0)= n^2(n-4)^2, \quad\quad \mathcal{Q}(1)= -2^{12},\quad\quad
\mathcal{Q} \left( \frac{n+2}{n-2} \right) = \frac{2^8n^2(n-4)^2}{(n-2)^4} \,,
\end{equation*}
while a short computation gives
\begin{equation*}
\mathcal{Q}(p_n) = \mathcal{Q}\left( \frac{n+4}{n-4} \right) = - \frac{2^{15}\,n^2}{(n-4)^3} \,.
\end{equation*}
This forces $\mathcal{Q}(p)$ to have three real roots in the interval $(0,p_n)$, so the fourth root must also be real. To
find its exact location, we compute
\begin{equation}\label{lim2}
\lim_{p\to \pm\infty} \frac{\mathcal{Q}(p)}{p^4} = (n-4) \cdot (n^3-4n^2-128n+256)
\end{equation}
and we examine two cases.

\Case{1} When $4<n\leq 12$, the limit in \eqref{lim2} is negative.  Since $\mathcal{Q}(0)$ is positive by above,
the fourth root lies in $(-\infty,0)$, so $\mathcal{Q}(p)$ is negative for any $p\geq p_n$ whatsoever.

\Case{2} When $n\geq 13$, the limit in \eqref{lim2} is positive.  Since $\mathcal{Q}(p_n)$ is negative by above, the
fourth root $p_c$ lies in $(p_n,\infty)$, so $\mathcal{Q}(p)$ is negative on $[p_n,p_c)$ and non-negative on
$[p_c,\infty)$.
\end{proof}

\begin{lemma}\label{use}
Let $n>4$ and $p>\frac{n+4}{n-4}$.  Set $m=\frac{4}{p-1}$ and let $Q_4$ be the quartic in \eqref{Q}. Then
\begin{equation*}
\mathscr{R}(\mu) = Q_4(m-\mu) - Q_4(m)
\end{equation*}
has four real roots $\mu_1 < \mu_2 < \mu_3 =0 < \mu_4$.
\end{lemma}

\begin{proof}
As in the proof of Lemma \ref{poly}, we exploit the fact that $\mathscr{R}(\mu)$ is symmetric about
\begin{equation*}
\mu_* \equiv m- \frac{n-4}{2} = \frac{4}{p-1} - \frac{n-4}{2} <0.
\end{equation*}
It is clear that $\mu_3=0$ is a root of $\mathscr{R}(\mu)$.  Then $\mu_2= 2\mu_*<0$ must also be a root by symmetry. To
see that a positive root $\mu_4>m$ exists, we note that
\begin{equation*}
\lim_{\mu\to +\infty} \mathscr{R}(\mu) = +\infty, \quad\quad \mathscr{R}(m)= - Q_4(m) < 0
\end{equation*}
by \eqref{phia}.  Then $\mu_1= 2\mu_*-\mu_4= \mu_2-\mu_4 <\mu_2$ must also be a root by symmetry.
\end{proof}

Finally, we address the intersection properties of the solutions provided by Lemma \ref{ex} in the supercritical case
$p\geq p_c$. To this end, let us also introduce the function
\begin{equation}\label{sing}
\Phi(x) = Q_4(m)^{\frac{1}{p-1}} \cdot |x|^{-m}, \quad\quad m=\frac{4}{p-1}
\end{equation}
which is easily seen to be a singular solution of \eqref{ee}.

\begin{theorem}\label{sep}
Suppose that $n\geq 13$ and $p\geq p_c$, where $p_c$ is given by Theorem \ref{main}.  In other words, suppose that $n>4$
and $p>\frac{n+4}{n-4}$ and that \eqref{sc} holds.  If $\phi_\al,\phi_\b$ are any two of the solutions provided by Lemma
\ref{ex}, then

\begin{itemize}
\item[(a)]
the graph of $\phi_\al$ does not intersect the graph of the singular solution \eqref{sing};

\item[(b)]
the graph of $\phi_\al$ does not intersect the graph of $\phi_\b$, unless $\al=\b$.
\end{itemize}
\end{theorem}

\begin{proof}
To establish part (a), we have to show that
\begin{equation}\label{bel2}
|x|^m \,\phi_\al(x) < Q_4(m)^{\frac{1}{p-1}}
\end{equation}
for each $x\in\R^n$.  This amounts to a slight refinement of inequality \eqref{bel1} in Proposition \ref{st}, as we now
need the inequality to be strict. Let us then consider the function
\begin{equation*}
W(s) = r^m \phi_\al(r), \quad\quad s=\log r=\log |x|, \quad\quad m=\frac{4}{p-1} \,.
\end{equation*}
Inequality \eqref{bel1} in Proposition \ref{st} reads
\begin{equation}\label{21}
W(s) \leq Q_4(m)^{\frac{1}{p-1}}, \quad\quad s\in\R
\end{equation}
and we now have to show that this inequality is actually strict.  Since
\begin{equation*}
\lim_{s\to -\infty} W(s) = \lim_{r\to 0^+} r^m\phi_\al(r) = 0 < Q_4(m)^{\frac{1}{p-1}}
\end{equation*}
by \eqref{phia}, we do have strict inequality near $s=-\infty$.  Suppose equality holds at some point, and let $s_0$ be
the first such point.  Since $W(s)$ reaches its maximum at $s_0$, we then have
\begin{equation}\label{end2}
W(s_0) = Q_4(m)^{\frac{1}{p-1}}, \quad\quad W'(s_0)=0, \quad\quad W(s) < Q_4(m)^{\frac{1}{p-1}}
\end{equation}
for each $s<s_0$.  When it comes to the interval $(-\infty,s_0)$, we thus have
\begin{align*}
W(s)^p - Q_4(m)^{\frac{p}{p-1}} > pQ_4(m)\cdot \left( W(s)-Q_4(m)^{\frac{1}{p-1}} \right)
\end{align*}
by convexity.  This is the same inequality as \eqref{12}, except that the inequality is now strict.  In particular, the
argument that led us to \eqref{14} now leads us to a strict inequality
\begin{equation*}
(\d_s - \la_1)(\d_s - \la_2)(\d_s - \la_3)(\d_s - \la_4) \,Y(s) > 0, \quad\quad s<s_0
\end{equation*}
for some $\la_1,\la_2,\la_3<0<\la_4$.  Using the same argument as before, we get
\begin{equation*}
(\d_s - \la_3)(\d_s - \la_4) \,Y(s) > 0, \quad\quad s<s_0
\end{equation*}
because $\la_1,\la_2<0$.  Multiplying by $e^{-\la_3s}$ and integrating over $(-\infty,s_0)$, we then get
\begin{equation*}
e^{-\la_3s_0} \cdot [Y'(s_0) -\la_4 Y(s_0)] > 0
\end{equation*}
because $\la_3<0$ as well.  In view of the definition \eqref{Y} of $Y(s)$, this actually gives
\begin{equation*}
W'(s_0) >\la_4 \left( W(s_0) - Q_4(m)^{\frac{1}{p-1}} \right),
\end{equation*}
which is contrary to \eqref{end2}.  In particular, the inequality in \eqref{21} must be strict at all points and the
proof of part (a) is complete.

In order to prove part (b), we shall first show that
\begin{equation}\label{end3}
W'(s) > 0, \quad\quad s\in\R.
\end{equation}
Using Lemma \ref{cco} and the strict inequality in \eqref{21}, we find that
\begin{equation*}
\Bigl[ Q_4(m-\d_s) - Q_4(m) \Bigr] \,W(s) = W(s)^p - Q_4(m) W(s) < 0.
\end{equation*}
Invoking Lemma \ref{use}, we now factor the left hand side to get
\begin{equation*}
(\d_s - \mu_1)(\d_s - \mu_2)(\d_s - \mu_4) \,W'(s) < 0
\end{equation*}
for some $\mu_1 < \mu_2 < 0 < \mu_4$.  Once again, the last equation easily leads to
\begin{equation}\label{16}
(\d_s - \mu_4) \,W'(s) < 0
\end{equation}
because $\mu_1,\mu_2<0$. This makes $W'(s)-\mu_4W(s)$ decreasing for all $s$, so the limit
\begin{equation*}
\lim_{s\to +\infty} W'(s) - \mu_4W(s)
\end{equation*}
exists.  Since $W(s)$ tends to a finite limit as $s\to +\infty$ by \eqref{phia}, it easily follows that $W'(s)$ must
approach zero as $s\to +\infty$.  Since $\mu_4>0$ by above, this gives
\begin{equation*}
\lim_{s\to +\infty} e^{-\mu_4s} \,W'(s) = 0,
\end{equation*}
and then we may integrate \eqref{16} over $(s,+\infty)$ to deduce the desired \eqref{end3}.

Now, the inequality \eqref{end3} we just proved can also be written as
\begin{equation}\label{31}
0< W'(s) = r\cdot \d_r [r^m \phi_\al(r)] = r^m \,[m\phi_\al(r) + r\phi_\al'(r)]
\end{equation}
in view of our definition \eqref{W}.  On the other hand, a scaling argument shows that the solutions provided by Lemma
\ref{ex} are subject to the relation
\begin{equation}\label{32}
\phi_\al(r) = \al\phi_1(\al^{1/m} r).
\end{equation}
Differentiating \eqref{32} and using \eqref{31}, one now easily finds that $\d_\al \phi_\al(r)>0$ for all $\al,r>0$.  In
particular, the graphs of distinct solutions cannot really intersect, as needed.
\end{proof}

\end{document}